\documentclass[11pt,a4paper]{article}  

\usepackage[a4paper, top=1in, left=1in, right=1in, bottom=1in]{geometry}

\date{}
\usepackage{braket}
\usepackage{graphicx}
\usepackage{amsmath,amssymb,amsfonts}
\usepackage{hyperref}
\usepackage{bbold}
\usepackage{caption}
\usepackage{amsmath}  
\usepackage{amssymb}   
\usepackage{braket}
\usepackage{amsthm}  
\usepackage{lipsum}

\usepackage{amssymb}
\usepackage{changepage}
\usepackage{breqn}
\usepackage{amsmath}
\usepackage{amsfonts}
\usepackage[utf8]{inputenc}
\usepackage[T1]{fontenc}
\usepackage{amsmath}
\usepackage{amsfonts}
\usepackage{amssymb}
\usepackage{dsfont}
\usepackage{graphicx}
\usepackage[utf8]{inputenc}
\usepackage{amssymb}

\usepackage{changepage}

\usepackage{changepage}
\usepackage{breqn}
\usepackage{amsmath}
\usepackage{amsfonts}
\usepackage{tikz}
\usepackage[mathscr]{euscript}
\usepackage{yfonts}
\usepackage{tikz-cd}

\usetikzlibrary{positioning}

\usepackage{graphicx}
\usepackage[utf8]{inputenc}
\usepackage{amssymb}
\usepackage{amsthm}  
\usepackage{subcaption} 
\usepackage{url}
\usepackage{tikz}

\usepackage{graphicx}
\usepackage{bm}
\usepackage{latexsym}
\usepackage{amsmath}
\usepackage{bbm}
\usepackage{dsfont}

\usepackage{graphicx}
\usepackage{bm}
\usepackage{latexsym}
\usepackage{amsfonts}
\usepackage{amssymb}
\usepackage{pax}
\usepackage{wasysym}

\usepackage{amsfonts}
\usepackage{amssymb}

\usepackage{pax}
\usepackage{wasysym}

\title{\large Monotonicity of the Cheeger constant under Ricci flow on spheres}

\begin{document}
\maketitle
\date{}
Hollis Williams$^{1,2,*}$

[1]{Department of Physics, University of Warwick, Coventry CV4 7AL, United Kingdom}

[2]{\textit{Current address:} Theoretical Sciences Visiting Program, Okinawa Institute of Science and 
Technology Graduate University, Onna 904-0495, Japan}
\newline
\newline
Contact email: $^*$holliswilliams@hotmail.co.uk
\newline

	\begin{abstract}\noindent
We study the behavior of the Cheeger isoperimetric constant under the
Ricci flow on compact surfaces.  For metrics on a surface diffeomorphic
to $S^2$, we show that the Cheeger constant is non-decreasing along the
flow.  The proof uses evolution identities for parallel curves together
with a viscosity formulation of the evolution of $\log h$ which accommodates for the possible switching of minimizing regions.  We
also give examples of nontrivial Ricci flows on topological $2$-spheres
for which the Cheeger constant remains constant, demonstrating that strict monotonicity
is not expected.

	\end{abstract}

\noindent
\section{Introduction}

There is a rich history of the interaction between the mean curvature flow and geometric inequalities which relate quantities such as area and volume.  As an example, the theory of the one-dimensional mean curvature flow (usually known as the curve shortening flow) was developed in part using the isoperimetric inequality.  One can also go in the other direction and use the flow theory to derive inequalities (the theory of the curve shortening flow can be used to derive an isoperimetric inequality which generalizes inequalities of Alexandrov, Huber and Bol) [1, 2].  Another flow which might be useful in this context is the Ricci flow.  Given a smooth closed manifold $M$ equipped with a smooth Riemannian metric $g$, the Ricci flow is the following PDE:

\[ \frac{\partial g}{\partial t} = - 2 \text{Ric} (g)                , \]

\noindent
where $  \text{Ric} (g) $ is the Ricci curvature of $g$.

Hamilton observed that certain isoperimetric quantities improve
along the Ricci flow on $S^2$, and in particular he proved monotonicity
of an isoperimetric ratio related to lengths of closed curves on the
round sphere [3, 4].  More refined
isoperimetric comparison results for the normalized flow on $S^2$
were later obtained by Andrews and Bryan, who showed
that the round sphere is optimal among all metrics with the same area [5].  In this article, we explore the behavior of the Cheeger isoperimetric constant $h(M,g)$ under Ricci flow, which is defined for a Riemannian surface by

\[
h(M,g) := \inf_{\gamma} 
  \frac{\operatorname{Length}_g(\gamma)}
       {\min\{\operatorname{Area}_g(M^+),\operatorname{Area}_g(M^-)\}},
\]
where the infimum is taken over smooth embedded curves which separate the
surface.  The Cheeger constant is a fundamental geometric quantity: it
controls the first eigenvalue of the Laplacian by the inequalities of Cheeger and Buser and plays an important role in
ergodic theory and the study of geodesic flows [6-8].  An open question asked by Manning in 2004 is whether the Cheeger constant improves under the Ricci flow on surfaces [9]. 


We give a partial answer in the case of
topological $2$--spheres: for a compact Riemannina surface diffeomorphic to $S^2$ and provided a natural algebraic condition which relates the minimizing length and the area, the Cheeger constant is non-decreasing under the Ricci flow.  The proof uses evolution identities for length and enclosed area along a
parallel foliation, together with a viscosity solution formulation
which accommodates for the fact that $h(M,g)$ is non-differentiable at times when the
minimizing side switches.  It is important that this result only partly answers Manning's question, since he was primarily interested in surfaces with negative curvature, whereas our arguments rely on the topology
of $S^2$.  For surfaces with
negative Euler characteristic, the corresponding evolution formulae lose the structure that they have in the spherical case and new ideas would be required to resolve the problem.  
Nevertheless, the case of spheres already shows that Ricci flow improves the Cheeger constant for a wide variety of examples.

In Section 3 we discuss examples showing that this monotonicity result
cannot generally be strengthened to include strict monotonicity.
Even amongst topological spheres there exist nontrivial Ricci flows whose
geometry evolves but for which the Cheeger constant remains constant moduluo a trivial rescaling of the unnormalized flow. We finish with conclusions and future directions in Section 4.

\section{The Cheeger constant under Ricci flow}

In this section, we will prove that the Cheeger constant is non-decreasing under the Ricci flow on a closed surface which is diffeomorphic to~$S^2$.  In Section 2.1, we  write down the basic definitions and in Section 2.2, we give the evolution identities which we will need for the Ricci flow.  For completeness, in Section 2.3 we provide derivations of the identities used in Section 2.2.  In Section 2.4, we give the evolution of $\log h$ when the minimizing side is fixed.  In Section 2.5, we show that  $\log h$ is a viscosity supersolution.  Finally, in Section 2.6 we give the precise statement of the result and explain how it follows directly from the discussion in the previous subsections.

\subsection{Definitions}
Let $(M,g)$ be a smooth, closed, oriented Riemannian surface.  If $\gamma\subset M$ is a smooth embedded closed curve that separates
$M$ into two open regions $M^+$ and $M^-$ (so $M\setminus\gamma=M^+\sqcup M^-$),
define
\[
h(\gamma;g)
  :=\frac{\operatorname{Length}_g(\gamma)}
          {\min\{\operatorname{Area}_g(M^+),\,\operatorname{Area}_g(M^-)\}}.
\]
The Cheeger constant of $(M,g)$ is defined to be
\[
h(M,g):=\inf_{\substack{\gamma\text{ smooth,}\\\text{embedded, separating}}}
h(\gamma;g).
\]
\newline
\newline
\noindent
In dimension two, the $(n-1)$-submanifolds in the usual definition are
just closed curves, so we adopt the above form for clarity.  By the work of Hass and Morgan, the infimum is attained:
for any closed surface, there exists a smooth embedded minimizing loop~$\beta$
with $h(\beta;g)=h(M,g)$ [10].  We use the notation
\[
A_\pm(\gamma;g)=\operatorname{Area}_g(M^\pm),\qquad
L(\gamma;g)=\operatorname{Length}_g(\gamma),
\]
and we drop the explicit dependence on $\gamma$ or $g$ when there is no ambiguity.

\subsection{Evolution identities under Ricci flow}
Let $g(t)$, $t\in I$, be a smooth solution of the Ricci flow on~$M$:
\[
\partial_t g(t)=-2\operatorname{Ric}(g(t))=-R(g(t))\,g(t),
\]

\noindent
The final equality follows because for a surface $\operatorname{Ric}=\tfrac12 R g$.
For any smoothly varying domain $\Omega(t)\subset M$ with smooth boundary
$\partial\Omega(t)$ we recall the standard first-variation formulas [3, 11]:

\begin{enumerate}
\item \textbf{Area of a fixed region:}
If $\Omega$ is a fixed region of the surface
\[
\frac{d}{dt}\operatorname{Area}_{g(t)}(\Omega)
   =-\!\int_{\Omega} R\,dA_{g(t)}.
\]
\item \textbf{Length of a fixed curve:}
If $\gamma$ is a fixed embedded curve,
\[
\frac{d}{dt}L_{g(t)}(\gamma)
   =-\!\int_{\gamma} \tfrac12 R\,ds_{g(t)}.
\]
\end{enumerate}

For a one-parameter family of parallel curves
$\{\gamma_\rho\}_{\rho\in(-\delta,\delta)}$
obtained by pushing $\gamma_0$ along a unit normal field, denote
$L(\rho,t)=L_{g(t)}(\gamma_\rho)$ and
$A_\pm(\rho,t)=\operatorname{Area}_{g(t)}(M^\pm_\rho)$.
The differential identities for the rate of change of these quantities with respect to time are  
\begin{align}
\partial_t\log L &=
   \partial_{\rho}^2(\log L)
   +\frac{\Gamma}{L}\,\partial_{\rho}(\log L)
   +Q_L(\rho,t), \label{eq:evolL}\\
\partial_t\log A_\pm &=
   \partial_{\rho}^2(\log A_\pm)
   +\frac{L^2}{A_\pm^2}
   -\frac{4\pi}{A_\pm}
   +\frac{\Gamma}{L}\,\partial_{\rho}(\log A_\pm)
   +Q_{A_\pm}(\rho,t), \label{eq:evolApm}\\
\frac{d}{dt}\log A &=
   -\frac{8\pi}{A}+Q_A(t), \label{eq:evolA}
\end{align}
where $A(t)=\operatorname{Area}_{g(t)}(M)$,
$\Gamma(\rho,t)=\int_{\gamma_\rho}k\,ds$ is the total geodesic curvature,
and the remainder terms $Q_L,Q_{A_\pm},Q_A$ vanish or are uniformly bounded
in a small tubular neighborhood where the foliation is smooth [11].
These identities isolate the principal second derivative and curvature terms
which control the sign in the evolution of the Cheeger constant $h$.

\subsection{Derivation of differential identities}

For completeness, we derive here the identities
\eqref{eq:evolL}–\eqref{eq:evolA}
used in Section 2 for the evolution of
the length of a parallel curve and the enclosed areas
under Ricci flow.  Let $\gamma_0\subset (M,g(t))$ be a smooth embedded closed curve
which separates $M$ into two components $M^+$ and $M^-$.  In a small tubular neighborhood of $\gamma_0$, we can introduce standard Fermi coordinates $(\rho,s)$,
where $\rho$ is the signed distance along the unit normal~$\nu$
to~$\gamma_0$, and $s$ is the arc length parameter along $\gamma_0$.  The metric takes the form
\[
g = d\rho^2 + f(\rho,s,t)^2\,ds^2,
\]
with $f(0,s,t)=1$ and $\partial_\rho f(0,s,t) = -k(s,t)$,
where $k$ denotes the geodesic curvature of $\gamma_0$.  The parallel curves are
$\gamma_\rho = \{ (\rho,s)\mid s\in[0,L_0) \}$,
and their length and enclosed areas are
\[
L(\rho,t) = \int_0^{L_0} f(\rho,s,t)\,ds, \qquad
A_+(\rho,t) = \int_0^\rho L(\sigma,t)\,d\sigma, \qquad
A(t) = A_+(t)+A_-(t).
\]

As stated above, the Ricci flow on a surface satisfies
$\partial_t g = -R\,g$,
so that
\[
\partial_t f = -\tfrac12 R\,f, \qquad
\partial_t(dA_g) = -R\,dA_g.
\]
We recall that the Gauss curvature satisfies the structure equation
\[
\partial_\rho k = -K - k^2,
\]
where $K$ denotes the Gaussian curvature of the surface and the scalar curvature is $R=2K$.  Differentiating the length of $\gamma_\rho$ in time gives
\begin{equation}\label{eq:dLdt}
\partial_t L(\rho,t)
   = \int_0^{L_0} \partial_t f(\rho,s,t)\,ds
   = -\tfrac12 \int_{\gamma_\rho} R\,ds,
\end{equation}
which is the standard first variation formula
used in Section 2.2.

Next, differentiate $L(\rho,t)$ with respect to $\rho$:
\[
\partial_\rho L(\rho,t)
   = \int_0^{L_0} \partial_\rho f(\rho,s,t)\,ds
   = -\int_{\gamma_\rho} k\,ds = -\Gamma(\rho,t),
\]
and hence
$\partial_\rho^2 L(\rho,t)
   = -\partial_\rho \Gamma(\rho,t)
   = \int_{\gamma_\rho} (K+k^2)\,ds$.
Dividing by $L(\rho,t)$ gives
\begin{equation}\label{eq:rhoL}
\partial_\rho^2 (\log L)
   + \frac{\Gamma}{L}\,\partial_\rho(\log L)
   = \frac{1}{L}\int_{\gamma_\rho}(K+k^2)\,ds.
\end{equation}

\noindent
Using $R=2K$ and combining
\eqref{eq:dLdt} and \eqref{eq:rhoL},
we eliminate the curvature integral to obtain
\[
\partial_t (\log L)
   = \partial_\rho^2(\log L)
     +\frac{\Gamma}{L}\,\partial_\rho(\log L)
     + Q_L(\rho,t),
\]
where $Q_L(\rho,t)$ collects
higher-order terms involving the variation of $k$ and $K$ along~$s$.  Note that these vanish if the metric and $\gamma_0$ are rotationally symmetric
or are uniformly bounded in a short collar.
This proves equation \eqref{eq:evolL}.

We next consider evolution of the enclosed area.  By definition,
$A_+(\rho,t)=\int_0^\rho L(\sigma,t)\,d\sigma$,
so differentiating twice with respect to~$\rho$ gives
\[
\partial_\rho^2 A_+ = \partial_\rho L
   = -\Gamma.
\]
Using the Gauss-Bonnet theorem for the domain
bounded by $\gamma_\rho$, we obtain
\[
\int_{M^+_\rho} K\,dA + \Gamma = 2\pi,
\]
and combining this with $\partial_t A_+ = -\int_{M^+_\rho}R\,dA
 = -2\int_{M^+_\rho} K\,dA$,
we get
\[
\partial_t A_+
   = -4\pi + 2\Gamma
   = -4\pi - 2\,\partial_\rho^2 A_+.
\]
After dividing by $A_+$ and rearranging in logarithmic form, we have 
\[
\partial_t \log A_+
   = \partial_\rho^2(\log A_+)
     +\frac{L^2}{A_+^2}
     -\frac{4\pi}{A_+}
     + Q_{A_+}(\rho,t),
\]
where $Q_{A_+}$ absorbs higher-order or non-uniform terms,
yielding equation \eqref{eq:evolApm}.
An analogous identity holds for the area $A_-$.  We use the convention that for a domain $M^+_\rho$ with smooth boundary
$\partial M^+_\rho=\gamma_\rho$ the geodesic curvature $k_g$ of
$\gamma_\rho$ is taken with respect to the outward unit normal
of $M^+_\rho$.  With this sign convention, Gauss-Bonnet for the
topological disc $M^+_\rho$ takes the form
\[
\int_{M^+_\rho} K\,dA + \int_{\gamma_\rho} k_g\,ds = 2\pi,
\]
and, recalling our notation $\Gamma(\rho,t)=\int_{\gamma_\rho}k_g\,ds$,
we have the above formula for the domain bounded by $\gamma_{\rho}$.

For the total area, we integrate $\partial_t A = -\int_M R\,dA = -8\pi$
on the $2$-sphere by Gauss–Bonnet to obtain
\[
\frac{d}{dt}\log A = -\frac{8\pi}{A} + Q_A(t),
\]
which is~\eqref{eq:evolA}.
The remainder term $Q_A$ vanishes vanishes when the total area is
measured on a single topological sphere and is non-vanishing only if one
rescales the metric or restricts to a subdomain.  Collecting the preceding calculations yields the evolution identities
\eqref{eq:evolL}–\eqref{eq:evolA} used in Section~2.
The terms $Q_L,Q_{A_\pm}$ and $Q_A$ may be controlled
by uniform curvature and geometry bounds in a small collar neighborhood and therefore do not affect the sign estimates
used in the proof of monotonicity.

\subsection{Evolution of \texorpdfstring{$\log h$}{log h} }

\noindent
We next consider the evolution of \texorpdfstring{$\log h$}{log h} when the minimizing side is fixed.  The term minimizing side refers to the region
$M^+$ (resp.\ $M^-$) for which $A_+(t)\le A_-(t)$, i.e.\ the side which actually realizes the minimum in the definition $h=L/\min\{A_+,A_-\}$.
When that side remains fixed in time, $h$ varies smoothly and the following
differential geometry calculation applies.  To deal with times where the sides switch, we will later appeal
to a standard viscosity solution argument.

Begin by fixing $t_0$ and suppose that at $t_0$ the Cheeger minimizer
$\beta=\gamma_{\rho=0}$ is achieved by the side~$M^+$,
i.e.~$A_+(0,t_0)\le A_-(0,t_0)$.
For $t$ near~$t_0$ this side remains minimizing, so
\[
h(t)=\frac{L(0,t)}{A_+(0,t)},\qquad
\log h = \log L - \log A_+.
\]
Subtracting~\eqref{eq:evolApm} from~\eqref{eq:evolL} gives
\begin{equation}\label{eq:logh_evol_fixed}
\partial_t\log h
   =\partial_{\rho}^2(\log h)
    +\frac{\Gamma}{L}\,\partial_{\rho}(\log h)
    +\Bigl(\frac{L^2}{A_+^2}-\frac{4\pi}{A_+}\Bigr)
    +R(\rho,t),
\end{equation}
where $R(\rho,t)$ denotes a linear combination of the uniformly bounded remainder terms.  The algebraic factor satisfies
\[
\frac{L^2}{A_+^2}-\frac{4\pi}{A_+}
   =\frac{1}{A_+}\Bigl(\frac{L^2}{A_+}-4\pi\Bigr)
   =\frac{1}{A_+}\Bigl(h\,\frac{L}{A_+}-4\pi\Bigr).
\]
Therefore, if
\begin{equation}\label{eq:h_condition}
h(t)<\frac{4\pi}{L(0,t)},
\end{equation}
then this term is strictly negative.
At the minimizing curve $\rho=0$ we have
$\partial_\rho\log h(0,t)=0$ and
$\partial_\rho^2\log h(0,t)\ge0$,
so evaluating~\eqref{eq:logh_evol_fixed} at $(0,t)$ yields
\[
\partial_t\log h(0,t)\ge0.
\]
It follows that whilst the minimizing side remains fixed,
$h(t)$ is non-decreasing.

Regarding stability of the minimizer, the infimum in the definition of the
Cheeger constant on a closed surface is attained by a smooth embedded
minimizing loop $\beta$ for each fixed metric.  If this minimizer is
nondegenerate (isolated as a critical point of the relevant area/length
functional), then standard arguments with the implicit function theorem imply
that the minimizer persists smoothly under small smooth perturbations
of the metric.  In particular, for a smooth Ricci flow, the same loop can
be followed for a short time interval until either (i) the minimizer
ceases to be unique or (ii) a switching time is reached.  The latter
possibility will be handled with viscosity solutions, hence the derivation in the fixed minimizer regime is
justified for all non-switching times.  At switching times, only
one-sided representations are available and the viscosity formalism is
used.

\subsection{Viscosity supersolution and the parabolic maximum principle}

\noindent
In this subsection, we will prove a lemma for viscosity supersolutions.  For clarity, we incorporate the algebraic zero--order term arising in
equation \eqref{eq:logh_evol_fixed} into the operator used in the viscosity
argument.  Define
\begin{equation}\label{eq:operator-full}
\mathcal{L}[u]
:=\partial_t u - \partial_{\rho}^2 u - \frac{\Gamma}{L}\partial_{\rho}u
   - \Bigl(\frac{L^2}{A_+^2}-\frac{4\pi}{A_+}\Bigr) \, u.
\end{equation}
With this convention, equation \eqref{eq:logh_evol_fixed} reads as
\[
\mathcal L[\log h] \ge R(\rho,t),
\]
where $R(\rho,t)$ collects the uniformly bounded remainder terms.  In the sequel, we therefore treat $\log h$ as a viscosity supersolution of
$\mathcal L[u]=0$ with a controlled source term.  When the algebraic
factor has a favorable sign, the
zero--order contribution is nonpositive and the comparison principle
applies.
\newline
\newline
\noindent
$\textbf{Lemma 2.1}$: Let $g(t)$ be a smooth Ricci flow on $M\simeq S^2$, and let $h(t)=h(M,g(t))$
be the Cheeger constant defined as in Section~2.1. Assume the evolution
identities \eqref{eq:evolL}--\eqref{eq:evolA} hold in a collar neighborhood
of a minimizing curve and that the remainder terms $Q_L,Q_{A_\pm},Q_A$ are
uniformly bounded there. Define the linear operator as in equation (8).  Then $\log h$ is a viscosity supersolution of $\mathcal{L}[u]=0$ in the following sense:
if $\phi$ is a smooth test function and $\log h-\phi$ attains a local minimum
at $(\rho_0,t_0)$ then $\mathcal{L}[\phi](\rho_0,t_0)\ge 0$.

\begin{proof}
This can be proved with a standard viscosity argument [11].  First, consider the case when $(\rho_0,t_0)$ is a point where the minimizing side
is locally fixed (so the same side realizes the $\min\{A_+,A_-\}$ in a
neighborhood of $(\rho_0,t_0)$). Locally, we may then write
$\log h=\log L-\log A_+$ (or $\log L-\log A_-$ if the minus side minimizes).
By the calculation leading to \eqref{eq:logh_evol_fixed} we have
\[
\mathcal{L}[\log h] \ge 0
\]
up to the uniformly bounded remainder terms, and hence at any classical touching point
by a smooth $\phi$ with $\phi(\rho_0,t_0)=\log h(\rho_0,t_0)$ and
$\phi\le\log h$ in a neighborhood we obtain $\mathcal{L}[\phi](\rho_0,t_0)\ge 0$.

Next, consider the switching case where $A_+=A_-$ at $(\rho_0,t_0)$. There are now two one-sided representations: $\log h=\log L-\log A_+$ and
$\log h=\log L-\log A_-$. If $\phi$ touches $\log h$ from below at
$(\rho_0,t_0)$, then $\phi$ also touches at least one of the one-sided functions
from below (otherwise, one could locally lower $\phi$ on one side and contradict
the minimality). Applying the one-sided computation to that representation
and again using the uniform control of the remainders gives
$\mathcal{L}[\phi](\rho_0,t_0)\ge0$.  This establishes the viscosity supersolution property.
\end{proof}
\noindent
$\textbf{Corollary 2.2}$: Under the hypotheses of Lemma 2.1, suppose in
addition that the sign condition
\[
h(t) < \frac{4\pi}{L(\beta,t)}
\]
holds at time $t$ for the minimizing curve $\beta$. Then $t\mapsto h(t)$ is
non-decreasing at that time. If the condition holds for all $t$ in an interval,
then $h(t)$ is non-decreasing on that interval.

\begin{proof}
By Lemma 2.1, $\log h$ is a viscosity supersolution
of $\mathcal{L}[u]=0$. The strict inequality implies that the zero-order
algebraic term in \eqref{eq:logh_evol_fixed} is nonnegative when evaluated
at the minimizing curve, so in the supersolution inequality we obtain
$\mathcal{L}[\log h]\ge 0$ in the viscosity sense with a nonnegative zero-order
contribution. The parabolic maximum principle for viscosity supersolutions implies that $\log h$ cannot
decrease in time; equivalently, $h(t)$ is non-decreasing [11, 12]. 
\end{proof}
\noindent
$\textbf{Proposition 2.3}$: Let $g(t)$ be a smooth solution of the Ricci flow on a closed surface
$M\simeq S^2$.
If at some time~$t_0$ the Cheeger minimizer~$\beta$
satisfies
\[
h(M,g(t_0))<\frac{4\pi}{L(\beta,g(t_0))},
\]
then $\partial_t h(M,g(t))|_{t=t_0}\ge0$.
Consequently, if this inequality holds for all times in an interval,
then $h(M,g(t))$ is non-decreasing on that interval.

\begin{proof}
Equation~\eqref{eq:logh_evol_fixed} and the variational conditions at the
minimizer imply $\partial_t\log h(0,t_0)\ge0$
under~\eqref{eq:h_condition}.
Nondifferentiability at switching times is handled by the viscosity
argument above, completing the proof.
\end{proof}

\subsection{Main theorem}
\noindent
Combining the local differential inequality with the viscosity supersolution property yields the following conditional monotonicity statement for the Cheeger constant.
\newline
\newline
\noindent
$\textbf{Proposition 2.4}$: Let $g(t)$ be a smooth solution to the Ricci flow on a closed surface
$M\simeq S^2$.
Fix a time $t_0$ and let $\beta$ be a smooth embedded curve achieving
$h(M,g(t_0))$.
Assume that in a tubular neighborhood of~$\beta$ the parallel foliation
$\{\gamma_\rho\}$ is smooth and that the curvature and geodesic curvature
are uniformly bounded so that the remainder terms
$Q_L,Q_{A_\pm},Q_A$ in
\eqref{eq:evolL}--\eqref{eq:evolA}
satisfy uniform bounds
\[
|Q_L|+|Q_{A_\pm}|+|Q_A|\le C_0
\]
for some constant $C_0$ independent of $\rho$ and $t$ in the collar.

If the inequality
\begin{equation}\label{eq:ineq4pi}
h(M,g(t)) < \frac{4\pi}{L(\beta,g(t))}
\end{equation}
holds at a given time $t$, then $\frac{d}{dt}h(M,g(t))\ge0$ at that time.
Consequently, whenever \eqref{eq:ineq4pi} remains valid along the Ricci flow,
the Cheeger constant $h(M,g(t))$ is non-decreasing in $t$.

\begin{proof}
We begin by noting that the curve $\beta$ is guaranteed to exist by results of Hass and Morgan [10].  The uniform bounds on $Q_L,Q_{A_\pm}$ and $Q_A$ ensure that the
leading-order evolution formula
\eqref{eq:logh_evol_fixed} holds with controlled error.
At the minimizing curve $\rho=0$, one has
$\partial_\rho\log h=0$ and $\partial_\rho^2\log h\ge0$
by the variational property of a minimizing loop.
Under the inequality~\eqref{eq:ineq4pi},
the zero-order curvature term
$\frac{L^2}{A_+^2}-\frac{4\pi}{A_+}$
is nonpositive,
so evaluating~\eqref{eq:logh_evol_fixed} at $(\rho,t_0)$ gives
$\partial_t\log h(0,t_0)\ge0$.
If the minimizing side changes at some later time,
the viscosity supersolution argument
(Lemma 2.1)
shows that $\log h$ satisfies $\mathcal{L}[\log h]\ge0$
in the viscosity sense.
By the parabolic maximum principle, this removes the possibility of any local decrease in~$h$,
which therefore remains non-decreasing.
\end{proof}
\noindent
The constant $4\pi/L$ in the inequality $h<4\pi/L$ originates from isolating the algebraic curvature term
in~\eqref{eq:logh_evol_fixed}.  Weaker geometric assumptions could replace it if one can control $\frac{L^2}{A_+^2}-\frac{4\pi}{A_+}$ directly.
In particular, Papasoglu’s bound
$h(M,g)\le C/\sqrt{\operatorname{Area}(M,g)}$
implies that \eqref{eq:ineq4pi} automatically holds
whenever the total area is sufficiently large
relative to the minimizing length~$L$ [13].  Papasoglu's estimate gives a constant $C>0$ such
that
\[
h(M,g)\le \frac{C}{\sqrt{A(M,g)}}.
\]
Combining this with the condition \(\,h<4\pi/L\,\) we obtain the
explicit sufficient criterion
\[
\frac{C}{\sqrt{A}} < \frac{4\pi}{L},
\]
hence the hypothesis $h<4\pi/L$ is satisfied whenever the total area exceeds $(C L/4\pi)^2$.

The uniform control on the terms
$Q_L,Q_{A_\pm}$ and $Q_A$ follows from bounded curvature and a positive
injectivity radius in a small tubular neighborhood of the minimizer.
These are standard consequences of smooth Ricci flow on compact surfaces.
Smoothness of the minimizing loop and the variational inequalities
$\partial_\rho\log h=0$ and $\partial_\rho^2\log h\ge0$
used in the proof are ensured by existence results of Hass and Morgan [11].  Since $g(t)$
 is smooth on a compact surface, curvature and geodesic curvature are bounded on a small tubular neighborhood of any smooth embedded loop; hence the remainder terms $Q_L$, $Q_{A_{\pm}}$ and $Q_A$
appearing are uniformly bounded. Moreover, the injectivity radius is positive on compact surfaces, so the parallel foliation exists for some uniform collar width [11].

To clarify further, the remainder terms represent higher-order contributions
coming from spatial derivatives of the metric coefficients in Fermi
coordinates.  Under a smooth Ricci flow on a compact surface, these
coefficients and all their spatial and time derivatives are uniformly
bounded on any compact time interval by standard theory.  Consequently, for any smooth
minimizing loop $\beta$ there exists a collar neighborhood
$U\ni\beta$ and $\varepsilon>0$ so that in Fermi coordinates the metric
coefficients (and finitely many derivatives) are uniformly bounded on
$U\times[t_0-\varepsilon,t_0+\varepsilon]$.  It follows that
\[
|Q_L|+|Q_{A_\pm}|+|Q_A|\le C_0
\]
on this collar for $t\in[t_0-\varepsilon,t_0+\varepsilon]$, for some constant $C_0$ depending only on the geometry of the collar and
the curvature bounds on the chosen time interval.  In particular, $C_0$
is independent of $\rho$ and $t$ within the chosen neighborhood.

\section{Counterexamples to strict monotonicity}

In Section 2 we proved that on surfaces diffeomorphic to $S^2$ the
Cheeger constant is non-decreasing along the Ricci flow, provided the
inequality
\[
h(M,g(t))<\frac{4\pi}{L(\beta,g(t))}
\]
holds for the minimizing curve $\beta$.  We may ask if the result could be strengthened to strict monotonicity.  In this section, we illustrate that
strict monotonicity fails in general, even within the topology of the
sphere, by presenting natural examples for which the Cheeger
constant is geometrically constant in time under Ricci flow.

\subsection{Ricci solitons}

A Ricci soliton is a solution of the Ricci flow of the form
\[
g(t)=\lambda(t)\,\Phi_t^* g_0,
\]
where $\Phi_t$ is a one-parameter family of diffeomorphisms and
$\lambda(t)$ is a scaling factor.  In such cases, the geometry evolves
non-trivially (unless $\Phi_t$ is the identity), but the flow is
self-similar.  Since the Cheeger constant is invariant under diffeomorphism
and scales homogeneously under scaling, we have
\[
h(M,\lambda g)=\lambda^{-1/2} h(M,g),
\]
In other words, the behavior of the Cheeger constant $h$ along a soliton is determined only by the factor
$\lambda(t)$.

On a compact surface diffeomorphic to $S^2$ the only gradient Ricci
soliton is the shrinking round sphere.  In this case
$\lambda(t)=(1-2t)$ up to normalization, and the length and area of any region
scale by $(1-2t)^{1/2}$ and $(1-2t)$, respectively.  As a result we have for any
separating loop $\gamma$,
\[
h(\gamma;g(t))
 = \frac{L(\gamma,g(t))}{A_\pm(\gamma,g(t))}
 = \frac{(1-2t)^{1/2} L(\gamma,g_0)}{(1-2t) A_\pm(\gamma,g_0)}
 = (1-2t)^{-1/2} \, h(\gamma;g_0),
\]
and minimizing in~$\gamma$ gives
\[
h(M,g(t)) = (1-2t)^{-1/2} h(M,g_0).
\]
Although the metric evolves non-trivially, the Cheeger constant
changes exactly by the scaling factor and therefore does not
exhibit any additional monotonic behavior beyond the trivial
self-similar rescaling. 


\subsection{Noncompact solitons}

Noncompact gradient solitons such as the cigar soliton or the Bryant soliton evolve non-trivially
under Ricci flow.  For examples like the cigar soliton, the Cheeger constant $h(M,g(t)) \equiv 0$ for all $t$, although in general, a complete noncompact surface may have a positive Cheeger constant.  These
examples do not fall under the compactness hypotheses of Section 2, but
they illustrate further that self-similar Ricci flow need not induce a
strict increase of isoperimetric quantities.

 \begin{figure}[tbp]
 \centering
 \includegraphics[width=0.2\linewidth]{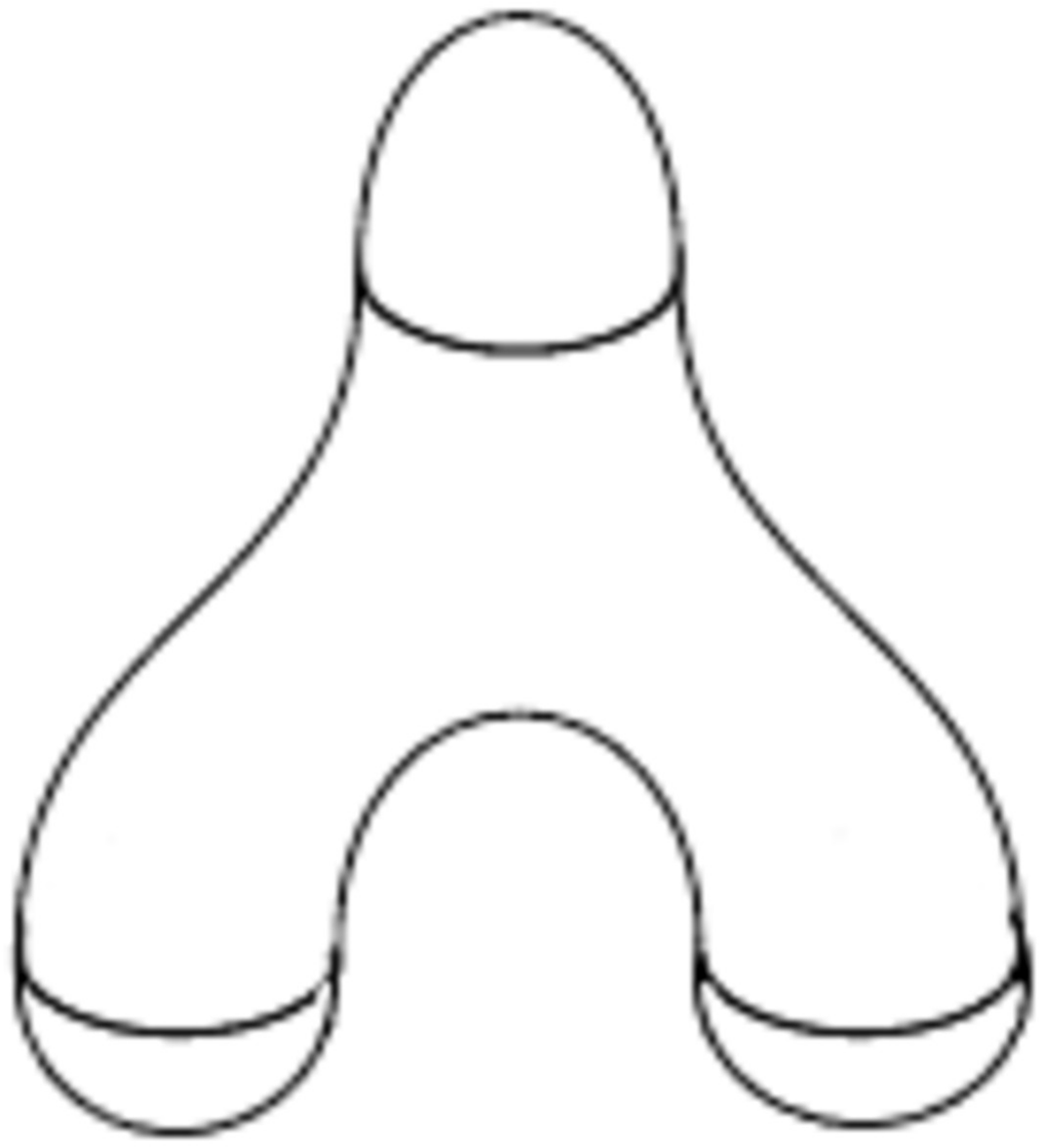}
 \caption{Pair of pants with spherical caps.  The caps have positive curvature, whereas the rest of the surface has negative curvature.}
 \label{fig:exmp}
 \end{figure}

\subsection{Pair of pants with spherical caps}

\noindent
We may also ask if there are examples of surfaces which have the topology of the sphere but exhibit regions of both positive and negative
curvature such that strict monotonicity can be broken.  One motivating model is a ``pair of pants'' surface whose
ends are capped with spherical pieces so as to obtain a smooth
topological 2-sphere (shown in Fig. 1).  Under Ricci flow such metrics evolve in a
non-self-similar manner and exhibit nontrivial curvature redistribution.
It is plausible in such geometries that the Cheeger minimizer may
undergo transitions for which $h(t)$ is constant on an interval or fails
to improve strictly.  A rigorous construction of this example would involve control of the evolving minimizing loops and the parallel foliation near the neck regions, which we leave for future work.

\noindent
\section{Conclusion}

In this article we proved that the Cheeger isoperimetric constant is
non-decreasing under the Ricci flow on any smooth compact surface which is diffeomorphic to the $2$--sphere.  The main new ingredient was a
viscosity formulation of the evolution of $\log h$, which
allows the classical argument for the isoperimetric ratio to remain valid when the minimizing side changes or
when the Cheeger minimizer is not smoothly differentiable in time.
Combined with the geometric evolution identities established in
Section~2 and the uniform control of the remainder terms in a tubular
neighbourhood of a minimizing loop, this yields a monotonicity
statement under the condition that
\[
h(M,g(t)) < \frac{4\pi}{L(\beta,g(t))}.
\]
As discussed, this inequality is generally satisfied whenever the total
area is sufficiently large relative to the minimizing length (for
example, by Papasoglu’s inequality for the Cheeger constant in terms of
area).

Section 3 shows that strict monotonicity cannot be expected, so that the above result is optimal.  In fact, even on topological $2$-spheres there exist nontrivial Ricci
flows whose geometry evolves but for which the Cheeger constant does not
increase except by the trivial rescaling of the metric.  In
particular, the shrinking round sphere evolves nontrivially under
unnormalised Ricci flow, but after renormalizing the metric to unit area, the Cheeger constant remains exactly constant in time.  No
geometric improvement of $h$ occurs along the flow, and strict monotonicity fails in
this natural sense.  Our result does not address the case of surfaces with negative Euler
characteristic studied by Manning and Katok [14]. Whether the
Cheeger constant should be strictly increasing under Ricci flow on
compact surfaces with negative curvature remains an open question.  The
methods developed here apply in the spherical setting because the
Gauss--Bonnet identity produces favorable sign conditions in the key
evolution formulas.  For surfaces of higher genus the corresponding
expressions have different signs, and additional ideas would be needed.  In summary, the analysis here shows that the Cheeger constant behaves in
a natural, geometrically controlled manner under Ricci flow on
spherical surfaces, providing a partial answer to the question of
how it might behave on a wider class of surfaces.

\section*{Acknowledgments}

\noindent
The author thanks Alexander Baumgartner for useful discussions.  This research was partly conducted whilst the author was visiting the Okinawa Institute of Science and 
Technology (OIST) through the Theoretical Sciences Visiting Program (TSVP).

\section*{References}

\noindent
[1]  P. Topping, Mean Curvature Flow and Geometric Inequalities, J. reine angew. Math. \textbf{503}, 47–61 (1998)

\noindent
[2] P. Topping, The isoperimetric inequality on a surface, Manuscripta Math. \textbf{100}, 23-33 (1999).

\noindent
[3] R. Hamilton, The Ricci flow on surfaces, Contemp. Math. \textbf{71}, 237-262 (1988).

\noindent
[4] R. Hamilton, An isoperimetric estimate for the Ricci flow on the two-sphere, Mod. Meth. in Comp. Ana. (Annals of Mathematics Studies, \textbf{137}, 191-200) (Princeton University Press, 1995).

\noindent
[5] B. Andrews and P. Bryan, Curvature bounds by isoperimetric comparison for normalized Ricci flow on the two-sphere, Calc. Var. \textbf{39}, 419-428 (2010). 

\noindent
[6] J. Cheeger, A lower bound for the smallest eigenvalue of the Laplacian, Problems in Analysis (Symposium in Honor of Salomon Bochner), 195-199 (Princeton University Press, 1970). 

\noindent
[7] P. Buser, A note on the isoperimetric constant, Ann. sci. de l’École Norm. Sup., 4 série, \textbf{15}(2), 213–230 (1982).

\noindent
[8] A. Katok, Lyapunov exponents, entropy and periodic orbits for diffeomorphisms, Pub. Math. de l’IHÉS \textbf{51}, 137–173 (1980).

\noindent
[9] A. Manning, The volume entropy of a surface decreases along the Ricci flow, Ergod. Theor. Dyn. Sys \textbf{24}, 171-176 (2004).

\noindent
[10] J. Hass and F. Morgan, Geodesics and soap bubbles in surfaces, Math. Zeit. \textbf{223}, 185-196 (1996).

\noindent
[11] B. Chow and D. Knopf, The Ricci flow: an introduction (American Mathematical Society, 2004).

\noindent
[12] M.G. Crandall, H. Ishii and P.-L. Lions, User's guide to viscosity solutions of second order partial differential equations, Bull. Amer. Math. Soc. (N.S.) \textbf{27}, 1-67 (1992). 

\noindent
[13] P. Papasoglu, Cheeger constants of surfaces and isoperimetric inequalities, Trans. Am. Math. Soc. \textbf{361}, 5139-5162 (2009).

\noindent
[14] A. Katok, Four applications of conformal equivalence to geometry and dynamics, Ergod. Th. Dynam. Sys. \textbf{8}, 139-152 (1988).


\end{document}